\documentclass[12pt]{article}
\usepackage{amsmath,amssymb}
\usepackage{color,url}
\usepackage{graphicx}    % standard LaTeX graphics tool

\def\blue#1{{\color{black}#1}}

\newtheorem{theorem}{Theorem}% [section]

\newtheorem{assumption}{Assumption}

\def\re{{\mathbb{R}}}

\def\cond{\hbox{\rm cond}}
\def\ocond{\overline{\hbox{\rm con}}\hbox{\rm d}}

\def\bee{\begin{equation}}
\def\ene{\end{equation}}
\def\beea{\begin{eqnarray}}
\def\enea{\end{eqnarray}}
\def\beeas{\begin{eqnarray*}}
\def\eneas{\end{eqnarray*}}

\def\mr#1{(\ref{#1})}

\def\ignore#1{}

\title{{\bf
    \blue{Iteration Complexity of Fixed-Step Methods by Nesterov and Polyak for
      Convex Quadratic Functions}}}

\author{
 \small Melinda Hagedorn, Heinrich Heine Univ., D\"usseldorf, Germany, melinda.hagedorn@hhu.de\\
 \small Florian Jarre, Heinrich Heine Univ., D\"usseldorf, Germany, jarre@hhu.de
}
\date  {Dec. 13, 2022,\medskip\\
Dedicated to Kees Roos and Florian Potra\medskip\\
All data generated or analysed in this article are available in \cite{hagedorn2}}
\begin{document}
\maketitle

\begin{abstract}
  This note considers the momentum method \blue{by Polyak and the
  accelerated gradient method by Nesterov, both}  without line search
  but with fixed step length applied to strictly convex quadratic functions
  assuming that exact gradients are
  used and appropriate upper and lower bounds for the extreme eigenvalues 
  of the Hessian matrix are known.
  Simple 2-d-examples show that the Euclidean distance of the iterates
  to the optimal solution is non-monotone. In this context
  an explicit bound is derived on the number of iterations needed to 
  guarantee a reduction of the Euclidean distance to the optimal solution
  by a factor $\epsilon$. \blue{For both methods
  the bound is optimal up to a constant factor, it complements earlier
  asymptotically optimal results for the momentum method,
  and it establishes another link of the momentum method and Nesterov's
  accelerated gradient method}.
\end{abstract}

\noindent
{\bf Key words:} Momentum method, convex quadratic optimization, iteration complexity

\section{Introduction}
\label{sec1}

Let $f:\re^n\to\re$ be a strictly convex quadratic function.
For the minimization of $f$ starting at an initial point $x^0\in\re^n$
with negative gradient $m^0:=-\nabla f(x^0)$
the following ``Momentum Method'' is considered for $k \ge 0$:
$$
x^{k+1} = x^k+\alpha m^k, \qquad m^{k+1} = \beta m^k - \nabla f(x^{k+1}).
\leqno(MM)
$$
Here $\alpha > 0$ is a sufficiently small step length that will be analyzed
below and $\beta\in[0,1)$ is a parameter that determines how much of the
previous search direction will be added to the new search direction. This 
parameter is also discussed below.

As pointed out for example in \cite{gower1}, with the initialization
$x^{-1}:=x^0$ the Momentum Method $(MM)$ 
is equivalent to the ``Heavy Ball Method'' (\cite{polyak}):
$$
x^{k+1} = x^k - \alpha\nabla f(x^{k}) +\beta(x^k-x^{k-1}).
\leqno(HBM)
$$ 

In spite of the proof of asymptotic optimality of the momentum method
$(MM)$ established almost 60 years ago in \cite{polyak} this method has
not received due attention until recently in the context of
machine learning and neural networks.
%, where it has
%found to be suitable also in the presence of noisy evaluations
%of the gradients.
It has proved to be very efficient in
numerical implementations, also when the gradients are replaced
with approximate stochastic gradients, see e.g. \cite{defazio_mm,gower1}.
Today modifications of the $(MM)$ are widely
used, in particular also in machine learning libraries.

\begin{figure}[h]
\begin{center}
  \caption{Convergence of the $(MM)$ for a 2-d-example with $\beta = 0.85$, and $\alpha = 1.9$.}
\label{pic.tp0}
%\vskip -10pt
%\begin{center}
%\phantom{.} \ \hskip -50pt
\includegraphics[width=170truemm]{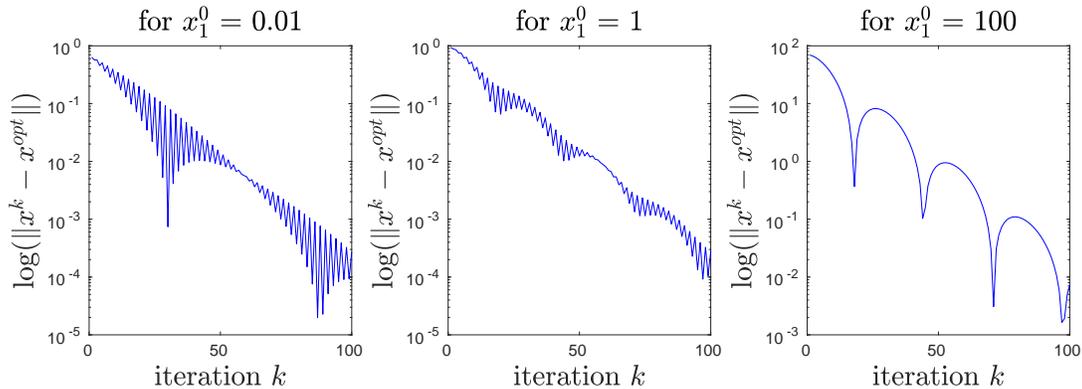}
\end{center}
%\vskip -285pt
\end{figure}

While being efficient in practice, the convergence behavior of momentum
methods is ``somewhat irregular''.
{\bf Figure \ref{pic.tp0}}
shows a typical non-monotone convergence behavior of the $(MM)$
for the simple 2-d-example  $f(x)\equiv\tfrac 12 (x_1^2+100x_2^2)$ with
$\beta = 0.85$ and \blue{$\alpha = \tfrac{1.9}{100}$}
during the first 100 iterations for the initial values $x^0=(\tfrac 1{100},1)^\text{T}$,
$x^0=(1,1)^\text{T}$, and $x^0=(100,1)^\text{T}$.

For convex quadratic functions the iterative process $(MM)$ or $(HBM)$
can be written as a recursion with a fixed matrix $M$, and the
spectral radius of $M$ is the main factor determining the
convergence behavior. A secondary factor is the approximation of
the spectral radius by a matrix norm.
The connection of both factors and the implication on the
choice of parameters is addressed in this paper.

\bigskip

\blue{
Closely related to the $(MM)$ is Nesterov's accelerated gradient method
\cite{nestg}. Following the presentation in \cite{nestg2} (Theorem 2.2.3),
it can be written as follows: Given $x^{0}\in\re^n$ set
$y^{0}:=x^{0}$ and for $k\ge 0$:
$$
x^{k+1}:=y^{k+1}+\beta (y^{k+1}-y^{k}) \qquad \hbox{where} \qquad
y^{k+1}:=x^{k}-\alpha \nabla f(x^{k})
\leqno(NAG)
$$
for suitable parameters $\alpha,\beta >0$.
As observed in \cite{yang}, for example, this method falls in the same
general framework as $(MM)$: Indeed, by eliminating the variable $y^k$
and setting $x^{-1}:=x^0$, the iteration
$(NAG)$ can also be written in the following compact form,
\bee\label{nag}
x^{k+1}:=x^{k}-\alpha \nabla f(x^{k}) +
\beta\left(x^{k}-x^{k-1} - \alpha(\nabla f(x^{k})-\nabla f(x^{k-1}))
\right),
\ene
where in contrast to $(HBM)$ the momentum term not only includes the
previous iterates $x^k$ but also
the previous descent steps ``$-\alpha\nabla f(x^{k})$''.
(To obtain the identical initilization as in $(NAG)$, in the very first step
$\nabla f(x^{-1})$ needs to be replaced with $0$.)
}

%\newpage
\subsection{Main Results}
The main theoretical
result of this paper is summarized in the following theorem.

\begin{theorem}\label{thm1}
  Let $f:\re^n\to\re$ be a strictly convex quadratic function
  with minimizer $x^*$ and
  denote $H:=\nabla^2f(x^*)$. Let
  $\overline{m}$ be an upper bound for the eigenvalues of $H$
  and let $0<\underline{m}$ be a lower bound for the eigenvalues
  of $H$. Then $\ocond(H):=\overline{m}/\underline{m}$
  is an upper bound for the condition number of $H$. Assume
  that $\ocond(H)\ge 28$.

  In $(HBM)$ define $\alpha=2/\overline{m}$, \ 
  \blue{$\beta = \left(1-\sqrt{2\ /\  \ocond(H)}\right)^2$
  and let  
  %The assumption on $\ocond(H)$ implies that $\delta \le e^{-2}$.
  $\epsilon\le 1\, / \, \ocond(H)$ be given.
  Then, given $x^0\in\re^n$, after
  $$
  1+\lceil \sqrt{2\ \ocond(H)}\
  \ln(\tfrac {2}\epsilon)\rceil
  $$
  steps of the $(HBM)$, an approximate solution
  $\bar x^k:=\tfrac 12 (x^{k-1}+x^k)$ is
  generated with
  $$
  \|\bar x^k-x^*\|_2\le
  %\sqrt{\|x^k-x^*\|_2^2+ \|x^{k-1}-x^*\|_2^2} \le
  \epsilon \|x^0-x^*\|_2.
  $$}
\end{theorem}

A corollary of the analysis of the $(MM)$ is the following
well known observation (see \cite{goh}, for example):
It may seem intuitive that with the use of a momentum $\beta >0$,
the step length $\alpha$ in $(MM)$ needs to be chosen more carefully.
However, at least for the minimization of
convex quadratic functions and in the absence of noise, the opposite is true:
 “If the step itself is made longer by adding more of the previous search
step to the new step (i.e.\ by increasing $\beta\in[0,1)$)
  then also the step length can be made longer
while maintaining global convergence.”

\bigskip

\blue{When $\overline{m}$ and $\underline{m}$ are known, then
  -- as detailed in Section~\ref{sec.nest} -- the analysis for $(MM)$
  also applies to $(NAG)$, and for an appropriate choice of
  $\alpha,\beta$, 
  an iteration complexity can be established for $(NAG)$
  that is a factor $\sqrt 2$ higher than the bound for $(MM)$
  in Theorem~\ref{thm1}. The factor $\sqrt 2$ arises since the step
  length is reduced by a factor 2 -- which essentially amounts to an
  increase of the estimate $\ocond(H)$ by a factor 2.}

\subsection{Related work}
The paper \cite{polyak} establishes for $(HBM)$ applied to
(non-quadratic smooth) strictly
convex functions local convergence of the form 
$$
\|x-x^*\|_2\le c(\epsilon) (1-\blue{\tfrac 2{1+\sqrt{\ocond(H)}}}+\epsilon)^k
\sqrt{\|x^1-x^*\|_2^2+ \|x^0-x^*\|_2^2}
$$
where $c(\epsilon)$ is a constant that only depends on $\epsilon$ 
(Theorem 9, Statement (3)). This result is obtained by analyzing the
spectral condition of a recursion associated with $(HBM)$ 
and using the fact that for a given matrix there
always is a vector norm and an associated matrix norm that
approximates the spectral radius to a given precision $\epsilon >0$.

For the case of convex quadratic functions the constant $c(\epsilon)$
is worked out explicitly in a slightly different setting in \cite{goh}.

The step length in Theorem~\ref{thm1} is limited such that
the first step is guaranteed not to increase the function value.
For larger values of $\ocond(H)$ the optimal step length $\alpha$
in \cite{polyak,goh} is nearly twice as long 
and leads to a substantial increase of the error components
associated with large eigenvalues of $H$ during the first several
iterations.  (This increase can
be reduced by taking at most half the step length of \cite{polyak,goh}
in the first iteration, a strategy that also applies to
methods with restart.)
As the generalization to non-quadratic functions also calls for
a shorter step length (\cite{gitman}, Footnote 3) and since
the use of long steps leads to an amplification of the
stochastic error associated with an inexact evaluation of the gradient,
subsequently, only step lengths at most $2/\overline m$ are considered
in this paper.
This imposes a restriction on $\alpha$ cutting off the
optimal step lengths identified in \cite{polyak,goh} at the
expense of a small constant factor in the theoretical efficiency
of the method.

Further studies of $(MM)$ for smooth strongly convex optimization
where the strong convexity parameter is to be estimated are
presented in \cite{barre,donoghue}.

Momentum methods and modifications thereof are in the focus of numerous
further recent research projects, in particular also for the non-convex case
-- which is much more difficult to analyze than the convex quadratic
case considered in this paper.
The articles quoted below form a rather brief and thus incomplete glimpse 
on recent lines of research on momentum methods.

An analysis comparing the $(MM)$ for the choice $\beta =0$ (steepest descent)
and $\beta>0$ in the non-convex case 
and in the presence of noise is presented, for example, in \cite{defazio_mm}.
This paper also gives a theoretical explanation for the
observed efficiency of $(MM)$.
The $(MM)$ is closely related to ADAM or ADAGRAD for which a recent analysis
covering the non-convex case is given in \cite{defossez}.
In particular, a new tight dependency on a certain ``heavy ball
momentum decay rate'' (which is zero in the context considered in this paper)
is established in \cite{defossez}.
Another recent approach that also considers the non-convex situation
without line search is analyzed in \cite{gratton,gratton2,gratton3}.
In this situation convergence to a second-order minimizer 
-- along with estimates on the rate of convergence --
is established under mild conditions.
An analysis generalizing $(HBM)$ and Nesterov's method
to a broad class of momentum methods and giving a unified
convergence analysis is given in \cite{diakonikolas}.\\
A unified analysis of stochastic momentum methods
in the weakly convex case -- and in the non-convex case --
is given in \cite{yang}, and a further unified analysis
that considers ``Quasi-Hyperbolic Momentum Methods''can be found
in \cite{gitman}. 
Limitations that arise when generalizing $(MM)$ to
stochastic optimization are addressed in the recent paper \cite{ganesh}.

%A comprehensive survey that also includes applications from deep learning and
%modifications of momentum methods for can be found in \cite{goodf}.

\section{Analysis of $(HBM)$}
The analysis of $(HBM)$ and thus also of $(MM)$ is
carried out with the following steps:
\begin{enumerate}
\item
  Section \ref{sec.spec} 
  follows the approach in \cite{polyak,donoghue} and derives
  a linear recursion for the iterates $x^k$ and analyzes 
  the spectral radius of the underlying system matrix.
\item
  Based on this analysis suitable parameters $\alpha,\beta$
  for $(HBM)$ are identified in Section~\ref{sec.par}.
\item
  Then the condition number of the similarity transformations
  to diagonalize the system matrix is analyzed.
\item
  It is shown that for the above parameters $\alpha,\beta$
  this condition number is unbounded in general.
\item
  A transformation to Schur canonical form is analyzed ``instead''.
\item
  Based on this transformation the main theoretical result is proved.
\end{enumerate}
The derivations in 
Steps 3) to Step 5) appear to be new; the results of
these steps, however, are essentially the same as in \cite{goh}
while the bound on the number of iterations derived in Step 6) may
not have been derived explicitly before.

As the $(HBM)$ of Section \ref{sec1} is invariant with respect to a
shift of the variable
and with respect to orthogonal transformations, for the analysis
it is assumed without loss of generality
that
\bee\label{assu}
f(x)\equiv \tfrac 12\, x^TDx
\ene 
with a positive definite diagonal matrix $D$.
  (The eigenvalue decomposition of the Hessian of $f$ that leads to
  the simple reformulation \mr{assu} is used only for the analysis
  of the algorithm, but not for the algorithm itself.)
In this case, the plain steepest descent method with $\beta = 0$ (no momentum)
converges if, and only if, $\alpha \in (0,\ 2/\max_{1\le i\le n} \{D_{i,i}\})$.

The following slightly weaker assumption will be made:
\begin{assumption}
It is assumed throughout that $f$ is given by \mr{assu} and
$$
\beta\in[0,1)\qquad\hbox{and}\qquad \alpha \in (0,\bar \alpha]
  $$
where $\bar \alpha := 2/\max_{1\le i\le n} \{D_{i,i}\}$
\end{assumption}

\bigskip

The $(HBM)$ 
can then be written with the following recursion
\bee\label{recursion}
\left(\begin{matrix} x^k \\ x^{k+1} \end{matrix} \right) =
\left(\begin{matrix} 0 & I \\ -\beta I & (1+\beta) I-\alpha D \end{matrix} \right)
\left(\begin{matrix} x^{k-1} \\ x^{k} \end{matrix} \right).
\ene
This is a discrete linear dynamical system
$\hat z^{k+1}=\hat M\hat z^k$ with the variable
$$
\hat z^k := \left(\begin{matrix} x^{k-1} \\ x^{k} \end{matrix} \right)
\qquad \hbox{and} \qquad
\hat M:= \left(\begin{matrix} 0 & I \\ -\beta I & (1+\beta) I-\alpha D \end{matrix} \right).
$$
Recursion \mr{recursion}  is block-separable, i.e.\
for $i\not=j$ the variables $x^k_i, \ x^{k+1}_i$
do not depend on $x^\ell_j$ for any $\ell \le k+1$.
Thus, when setting
$$
z^k_{(i)}:=\left(\begin{matrix} x^{k-1}_i \\ x^{k}_i \end{matrix} \right)
$$
for $k\ge 0$ and $1\le i\le n$, then Recursion \mr{recursion} can be written as
$$
z^{k+1}_{(i)} = M^{(i)}z^k_{(i)}\qquad \hbox{for} \ 1\le i\le n
$$ 
with
\bee\label{defmi}
M^{(i)}:=\left(\begin{matrix} 0 && 1 \\ -\beta  && 1+\beta -\alpha D_{i,i} \end{matrix} \right).
\ene

(This means that  
rows and columns of $\hat M$ can be permuted so that a block-diagonal matrix $M$
is obtained with $2\times 2$ diagonal blocks $M^{(i)}$ on the diagonal.)

\subsection{The spectral radius of $M^{(i)}$ in
  dependence of $\alpha$ and $\beta$}\label{sec.spec}

Possible convergence of the iterates $\hat z^k$ depends on
the norm $\|(M^{(i)})^k\|_2$ for large $k$, i.e.\ on
the spectral radius of $\hat M$
which coincides with the maximum spectral radius $\rho(M^{(i)})$
of $M^{(i)}$ for all $i$.

In the following the index $i$ is kept fixed. To simplify the notation, let
$$
\alpha_i:=\alpha D_{i,i}\in\,(0,2],
\quad \beta_i:=1+\beta-\alpha_i\in [-1,2),
\quad  \hbox{and} \ \gamma_i:=\sqrt{\beta_i^2-4\beta}.
$$
Here, $\gamma_i$ is either a non-negative real number
or a (purely imaginary) number with positive imaginary part.
In both cases, $\gamma_i^2=\beta_i^2-4\beta$ will be used below.
With these abbreviations $M^{(i)}$ is given by
$M^{(i)}=\begin{pmatrix}0&1\\-\beta&\beta_i\end{pmatrix}$.
Let further 
\bee\label{eigenv}
 \lambda_+ 
 :=\tfrac12(\beta_i+\gamma_i), \quad
 \lambda_- 
 :=\tfrac12(\beta_i-\gamma_i), \quad 
 v_+:=\begin{pmatrix}1\\ \lambda_+\end{pmatrix}, \
 \text{ and }\  v_-:=\begin{pmatrix}1\\\lambda_-\end{pmatrix}.
\ene
Observing that
$$
\lambda_\pm^2 \ = \
\tfrac 14 \beta_i^2 \pm \tfrac 12 \beta_i\gamma_i +\tfrac 14 \gamma_i^2 \ = \
\tfrac 14 \beta_i^2 \pm \tfrac 12 \beta_i\gamma_i +\tfrac 14 \beta_i^2-\beta
\ = \ -\beta + \tfrac 12 \beta_i^2 \pm \tfrac 12 \beta_i\gamma_i
\ = \ -\beta + \beta_i\lambda_\pm
$$
it follows that $M^{(i)}v_\pm=\lambda_\pm v_\pm$.
Thus, the -- possibly complex --
eigenvalues of $M^{(i)}$
are given by $\lambda_+$ and $\lambda_-$ and 
the spectral radius $\rho(M^{(i)})$ of $M^{(i)}$ is given by the larger
one of the two values
$$
|\lambda_\pm|\ =\ \tfrac12|\beta_i\pm \gamma_i| \ = \ 
\frac 12 \left| \beta_i
\pm\sqrt{\beta_i^2-4\beta} \right|.
$$
If the square root
is purely imaginary, i.e., if $\beta_i^2<4\beta$ then
the square of the absolute value is given by the sum of the squares of real and imaginary part,
i.e.
$$
\rho(M^{(i)})^2=\frac 14 \left(\beta_i^2 + |\gamma_i|^2)\right)=
\frac 14 \left(\beta_i^2 + (4\beta-\beta_i^2)\right)=\beta.
$$
Thus, the expression for $\rho(M^{(i)})$ simplifies to
\bee\label{defrho}
\rho:=\rho(M^{(i)}) = \left\{\begin{array}{ll}
\sqrt{\beta} & \hbox{if} \ \ \beta_i^2<4\beta\\
\frac 12\left( | \beta_i |+\gamma_i \right)
& \hbox{else}.
\end{array}\right.
\ene

Note that there is a double eigenvalue $\lambda_+=\lambda_-$
if $0=\gamma_i=\beta_i^2-4\beta=(1+\beta -\alpha D_{i,i})^2-4\beta$
with $\beta\in[0,1)$, i.e.\ if,
and only if, $\beta = (1-\sqrt{\alpha D_{i,i}})^2$ where
the definition of $\bar \alpha$ in Assumption~\ref{assu}
implies that $\alpha D_{i,i}\in(0,2]$ for all $i$.
(This is sketched on the right of Figure~\ref{specs} below.)

%\vfill\eject

\begin{figure}[h]
\begin{center}
  \caption{Spectral radius of $M^{(i)}$ as a function of $\alpha D_{i,i}\in (0,2]$
    and $\beta\in[0,1)$.}
\label{pic.tp1}
%\vskip -10pt
%\begin{center}
%\phantom{.} \ \hskip -50pt
\includegraphics[width=130truemm]{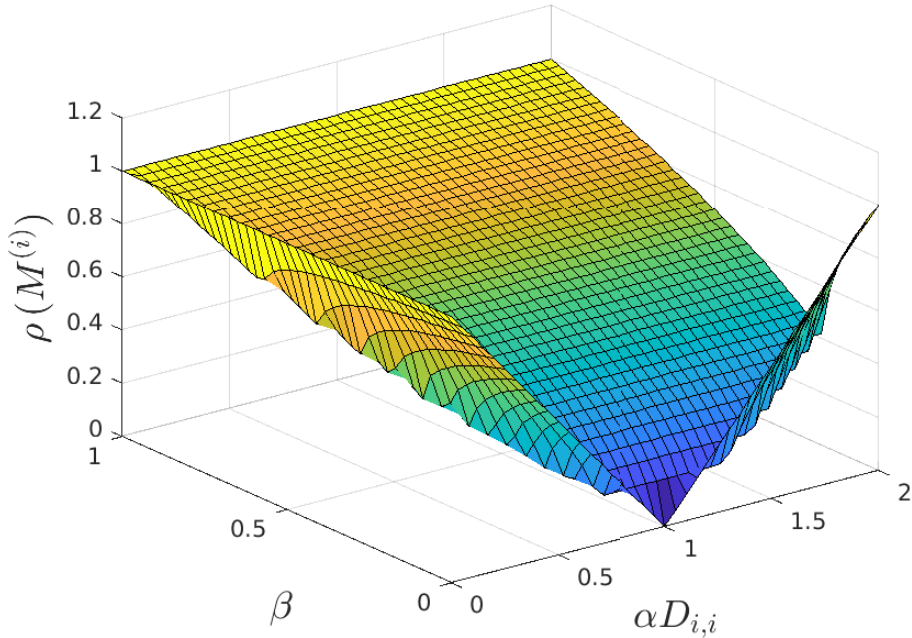}
\end{center}
%\vskip -285pt
\end{figure}

\bigskip

In {\bf Figure \ref{pic.tp1}}, 
the value of $\alpha D_{i,i}\in(0,2]$ is plotted from the middle to the right
and the value of $\beta$ from the middle to the rear-left. 
The associated values of $\rho$ are plotted on the vertical axis.

It is assumed that an upper bound $\overline{m}$
for the eigenvalues $D_{i,i}$ is known
and that $\alpha$ is chosen less than or equal to $2/\overline{m}$ so that 
the possible values of $\alpha D_{i,i}$  vary over a (typically wide)
range in the interval (0,2] and depend on the distribution of the eigenvalues.
This range is problem dependent and not subject to the design of the method.

A fast convergence of the HBM (or of the momentum method MM) is obtained if the value of
$\beta$ is chosen such that $\rho$ is small for a wide range of values
$\alpha D_{i,i}$. {\bf Figure \ref{pic.tp1}} gives the impression that this is achieved
when choosing $\beta = 0$. And indeed, when the values of
$\alpha D_{i,i}$ all cluster about the value ``1'', then the problem of
minimizing $f$ not only is quite easy (since the condition number of $D$
is close to 1) but also choosing $\beta = 0$ is nearly optimal.

\bigskip

\begin{figure}[h]
\begin{center}
  \caption{$\rho(M^{(i)})$ as a function of $\alpha D_{i,i}\in (0,2]$
    for fixed values of $\beta$.}
\label{pic.tp2}
%\vskip -10pt
%\begin{center}
%\phantom{.} \ \hskip -50pt
\includegraphics[width=130truemm]{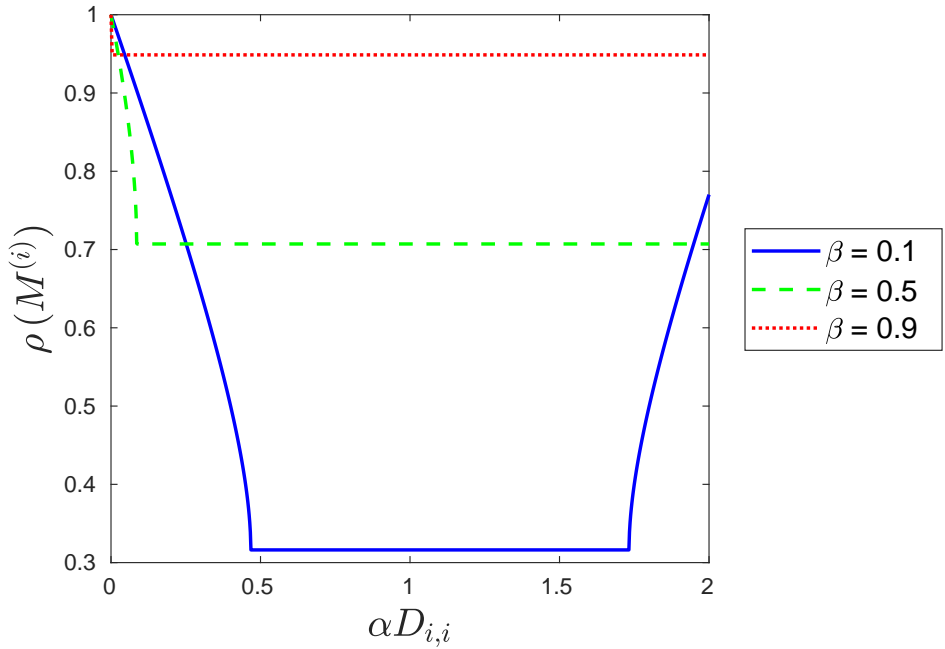}
\end{center}
%\vskip -285pt
\end{figure}

However, when the condition number of $D$ is poor, then very small
values of $\alpha D_{i,i}>0$ will occur. In this case, {Figure \ref{pic.tp1}}
is not suitable for understanding the best possible choice of
$\beta$. Instead, {\bf Figure \ref{pic.tp2}} displays 
the plots of $\rho(M^{(i)})$ as a function of $\alpha D_{i,i}\in (0,2]$
for $\beta \equiv 0.9$ in red, $\beta \equiv  0.5$ in green,
and $\beta \equiv  0.1$ in blue.

\begin{figure}[h]
\begin{center}
  \caption{$\rho(M^{(i)})$ as a function of $\alpha D_{i,i}\in (0,0.1)$
    for fixed values of $\beta$.}
\label{pic.tp3}
%\vskip -10pt
%\begin{center}
%\phantom{.} \ \hskip -50pt
\includegraphics[width=130truemm]{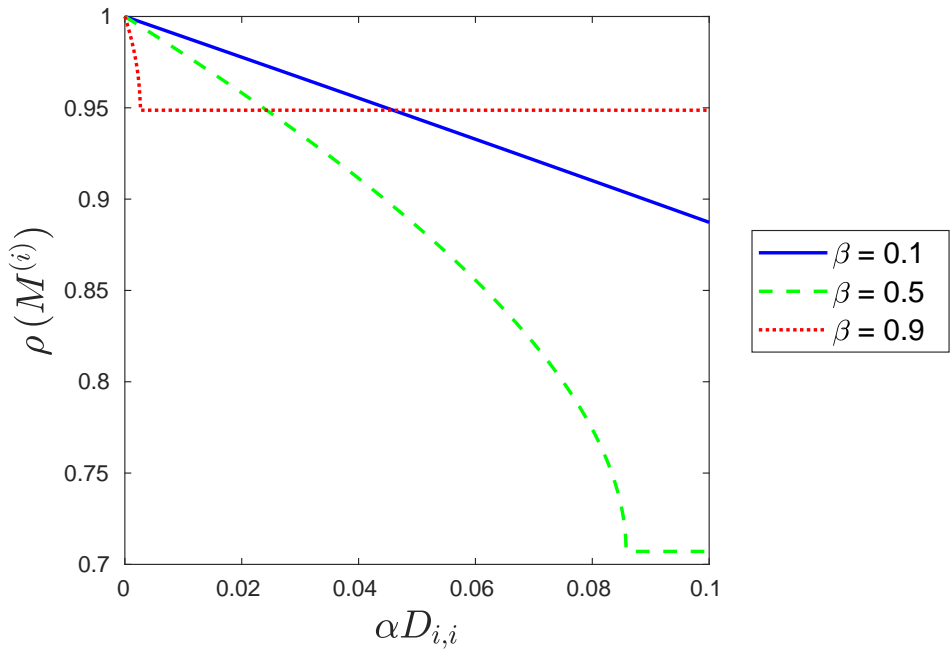}
\end{center}
%\vskip -285pt
\end{figure}

In {\bf Figure \ref{pic.tp2}} the choice $\beta=0.1$ in blue leads to
small values of $\rho(M^{(i)})$ when $\alpha D_{i,i}\in [0.5,1.5]$
(and a little beyond this interval).
For values $\alpha D_{i,i}\approx 0.2$, the choice $\beta=0.5$ in green
results in a much smaller value of  $\rho(M^{(i)})$ than $\beta=0.1$,
and as displayed in {\bf Figure~\ref{pic.tp3}}, when considering poorly
conditioned problems with some values of $\alpha D_{i,i}\approx 0.004$
the choice of $\beta=0.9$ in red results in a smaller value of
$\rho(M^{(i)})$ than $\beta=0.5$ or $\beta=0.1$.
({Figure \ref{pic.tp3}} is a zoom of {Figure \ref{pic.tp2}};
of course, small values of $\alpha D_{i,i}$ occurring in
ill-conditioned problems result in an increase of $\rho(M^{(i)})$
but for small $\alpha D_{i,i}>0$
the increase is much less for larger values of $\beta < 1$.)

  An interesting observation that can be deduced from
  {\bf Figure \ref{pic.tp2}} is the following:
  Recall that a step length $\alpha D_{i,i}>2$  is
  ``too long'' in the sense that it results in a divergent algorithm
  for the plain steepest descent method.
  For larger values of $\beta\ge 0.5$, however, the green line and the
  red line seem to (and do) continue ``for a little while'' with constant
  values to the right of $\alpha D_{i,i}=2$ which implies
  that for larger values of $\beta$ a step length that is
  a little bit ``too long'' still results in a convergent algorithm.
  In contrast to the intuition that with the use of momentum the step length
  needs to be chosen more carefully, actually the opposite is true.

\subsection{Selecting $\alpha$ and $\beta$}\label{sec.par}
In the following it is assumed that a lower bound
$$
0< \underline{m} \le \min_{1\le i\le n} D_{i,i} \qquad
\hbox{and an upper bound} \qquad \overline{m} \ge \max_{1\le i\le n} D_{i,i}
$$
are known.
Given $\overline{m}$, the step length $\alpha$ is fixed
to
\bee\label{alphapol}
\alpha := 2/\overline{m}.
\ene
In the following let
$$
\ocond(D)\ \ :=\ \ \overline{m}/\underline{m}\ \ \ge \ \ \cond(D)
$$
be an upper bound for $\cond(D)$.
Often, an upper bound for $\overline{m}$ can be obtained, for example by 
Gershgorin's theorem, while it may be difficult to determine
a positive lower bound for $\underline{m}$.
(The numerical effects of overestimating or underestimating
the condition number of $H$ are analyzed in \cite{donoghue}
along with an approach to adaptively adjust the estimate of the
condition number, provided that function evaluations are
available -- an assumption that unfortunately is not satisfied
in many stochastic settings.)

To simplify the discussion below, it is further assumed that
\bee\label{condd}
\ocond(D) > 2.
\ene
(In Theorem \ref{thm1} a stronger restriction on $\ocond(D)$
is made; for now \mr{condd} is sufficient.
For well-conditioned problems with $\ocond(D) \le 2$ the
plain steepest descent method with $\beta = 0$ is
optimal up to a factor of at most 2.)

%(This step length may be too long for the plain steepest
%descent method to converge, but as detailed below it is suitable
%for the method with momentum.)

When $\beta < 3-2\sqrt 2\approx 0.1716$ then there are
two values of $\alpha D_{i,i}$ for which $(1+\beta-\alpha D_{i,i})^2=4\beta$
(corresponding, for example, to the end points of the blue horizontal
line in Figure~\ref{pic.tp2}).
The considerations below apply to values
$$
\beta \ge 3-2\sqrt 2
$$
which have proved to be efficient in numerical experiments,
see, for example, \cite{gower1}, and
where only the left end point of the interval with constant values $\rho$ 
(i.e.\ of the horizontal lines in Figure~\ref{pic.tp2}) is relevant.

Let $\underline{\alpha} = \alpha\underline{m}$ be a lower bound for
possible values of $\alpha D_{i,i}$.
Assumption \mr{condd} implies 
$\underline{\alpha} < 1$.
Setting $\beta := (1-\sqrt{\underline{\alpha}})^2$, the spectral
radius of all $M^{(i)}$ (and thus, $\rho(M)$)) is upper
bounded by $\sqrt{\beta} = 1-\sqrt{\underline{\alpha}}$, i.e.
$$
\rho(M)\leq\sqrt{\beta} = 1-\sqrt{\underline{\alpha}} =
1-\frac {\sqrt 2}{\sqrt{\ocond(D)}}.
$$
When $\ocond(D)$ approximates $\cond(D)$ up to a
constant factor (this is what is meant with ``appropriate'' bound
in the abstract of this paper) then the resulting 
order $1- \frac 1{O(\sqrt{\cond(D)})}$ is the same order as 
the optimal rate of convergence of the conjugate gradient method,
see e.g. \cite{li}.
Moreover, the $2\times 2$ transformation matrices
transforming $M^{(i)}$ to diagonal form
or, more generally, to the Jordan canonical form
are independent
of the dimension, so that the 2-norm of the transformation matrices
that transform $\hat M$ to a canonical
form also is independent of $n$.\\
However, the optimal bound for the conjugate gradient iterations
applies to the
$D$-norm, $\|z\|_D=\sqrt{z^TDz}$, while the above rate applies to
a transformed problem, and the transformation matrix may be very
ill-conditioned. This is studied
further in Section \ref{sec.trans} below.

\subsection{On the norm of the transformation matrices}\label{sec.trans}

Consider the case in  \mr{eigenv} that  $\gamma_i=\sqrt{\beta_i^2-4\beta}$
is nonzero. Then $v_+\not=v_-$
and $M^{(i)}$ is diagonalizable, i.e.\ $M^{(i)}=S\Lambda S^{-1}$ where
$\Lambda$ is the diagonal matrix with the eigenvalues $\lambda_\pm=
\tfrac 12 (\beta_i\pm\gamma_i)$ of \mr{eigenv}
and 
$S:=S^{(i)}:=\begin{pmatrix}1&1\\ \lambda_+&\lambda_-\end{pmatrix}$
with columns $v_\pm$ as defined in \mr{eigenv}.

\bigskip

There are two cases: Either $\gamma_i>0$ is real or $\gamma_i$ is purely imaginary.

\bigskip

If $\gamma_i>0$ is real then $\lambda_+\lambda_- =
\tfrac 14 (\beta_i^2-\gamma_i^2)=\tfrac 14 (\beta_i^2-(\beta_i^2-4\beta))= \beta$ and
$$
S^TS = \begin{pmatrix}1+\tfrac 14 (\beta_i+\gamma_i)^2 &1+\beta \\
 1+ \beta  &  1+\tfrac 14 (\beta_i-\gamma_i)^2 \end{pmatrix}.
$$
Denote the eigenvalues of $S^TS$ by $\mu_+ \ge \mu_- > 0$.
Then 
\bee\label{mu1}
\mu_\pm=1+\tfrac 14 (\beta_i^2+\gamma_i^2)\pm
\tfrac 12 \sqrt{\beta_i^2\gamma_i^2+4(1+\beta)^2}.
\ene

\bigskip

If $\gamma_i$ is purely imaginary and $S^H$ denotes the complex conjugate
transpose of $S$, then 
$$
S^HS = \begin{pmatrix}
  1+|\lambda_+|^2 & 1 + \bar\lambda_+\lambda_-  \\[4pt]
  1 + \lambda_+\bar\lambda_-   & 1+|\lambda_-|^2
\end{pmatrix}.
$$
Using the definitions of $\lambda_\pm$ and $\gamma_i$ it follows
$\lambda_+\bar\lambda_-  =  \lambda_+^2  $ and $|\lambda_+|^2=|\lambda_-|^2=\beta$,
and $S^HS$ is given by
$S^HS=  \begin{pmatrix}
  1+\beta & 1 + \lambda_-^2  \\
  1 + \lambda_+^2   & 1+\beta
\end{pmatrix}$
with the eigenvalues
\bee\label{mu2}
\mu_\pm =1+\beta\pm |1+\lambda_+^2|.
\ene

\bigskip

\noindent
In both cases the condition number of $S$ is
$\cond(S)= \sqrt{\mu_+/\mu_-}$
with $\mu_\pm$ defined in \mr{mu1}~or~\mr{mu2}.

\bigskip

\begin{figure}[h]
%\begin{center}
  \caption{$\min\{\cond(S^{(i)}),\, 20\}$ 
    as a function of $\alpha D_{i,i}\in (0,2]$
    and $\beta\in[0,1)$ on the left, and zeros of the terms
      $\gamma_i$ depending on $\alpha D_{i,i}$
    and $\beta$ on the right plot.}
\label{specs}
\includegraphics[width=100truemm]{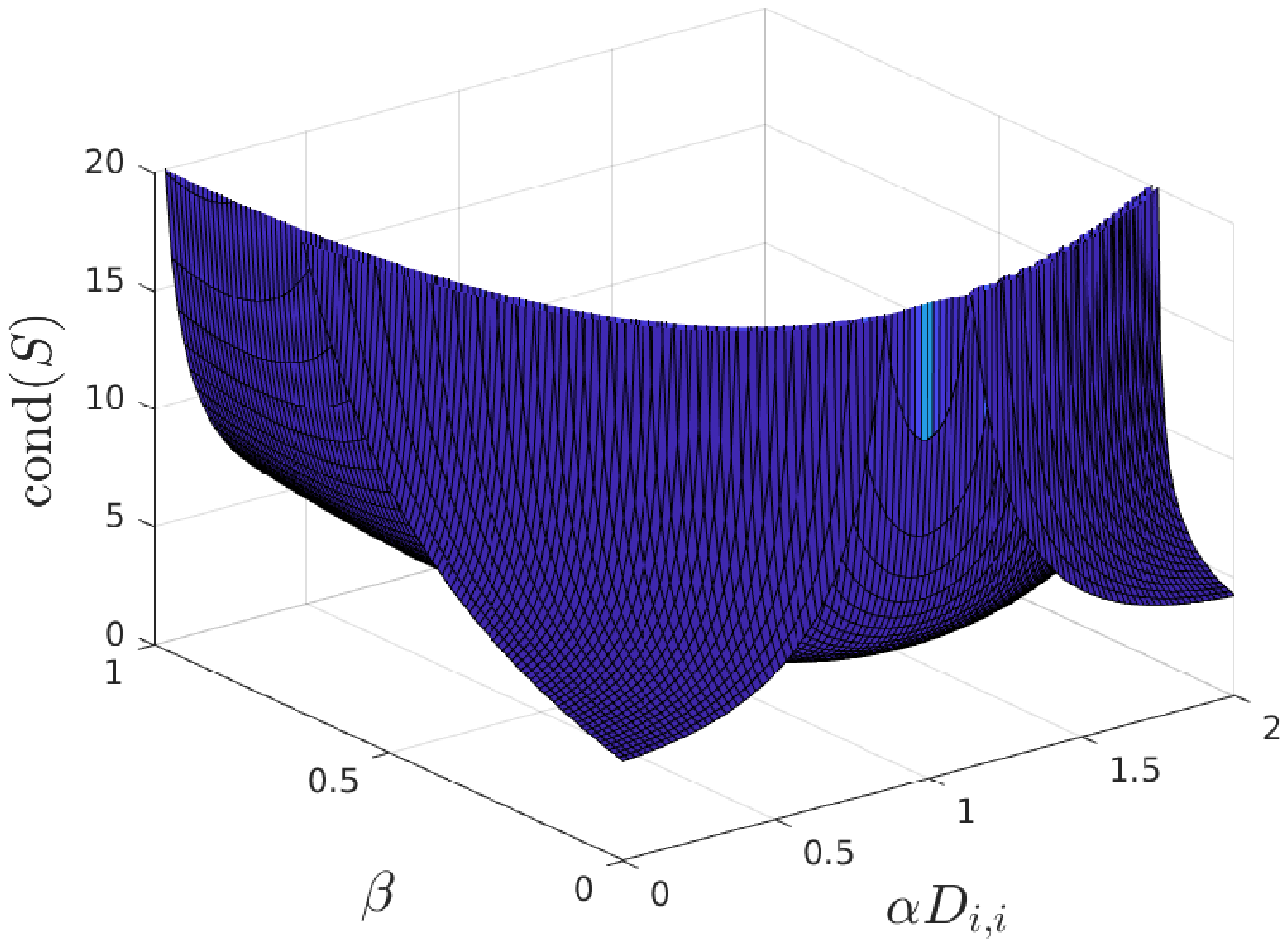}
\includegraphics[width=70truemm]{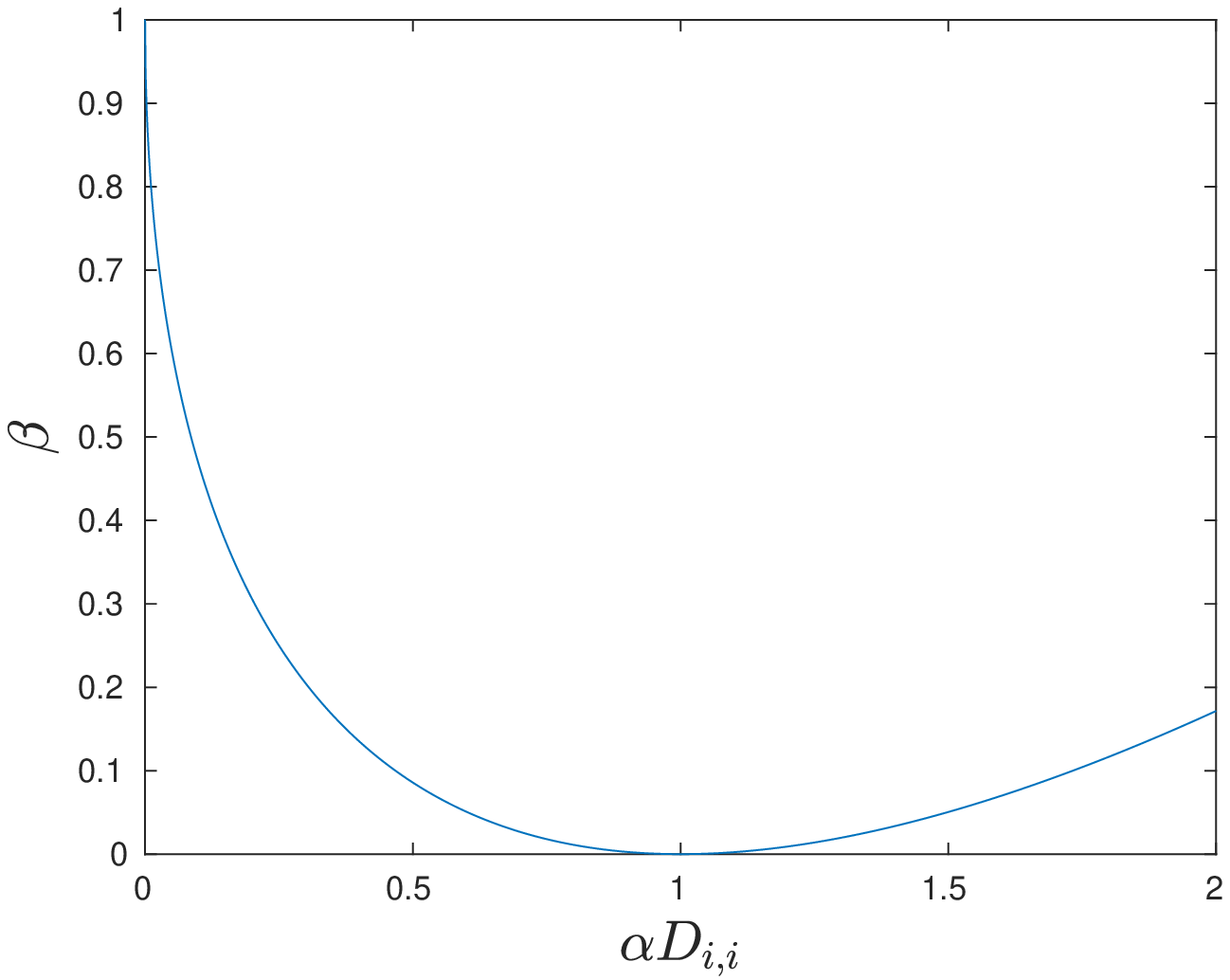}
%\end{center}
\end{figure}

Let $\Lambda$ denote the diagonal matrix with
diagonal entries $\lambda_+$ and $\lambda_-$.
As $\lambda_+$ and $\lambda_-$ approach each other,
i.e.\ as $\gamma_i \to 0$,
the columns of $S=S^{(i)}$ become nearly linearly dependent and 
$\cond(S)$ tends to infinity.
This observation is illustrated in {\bf Figure~\ref{specs}},
and it implies that the argument
\bee\label{argu1}
\|(M^{(i)})^k\|_2 = \|S\Lambda^kS^{-1}\|_2\le  \|S\|_2\|\Lambda^k\|_2\|S^{-1}\|_2
=\cond(S)\|\Lambda\|_2^k
\ene
for the convergence of $(M^{(i)})^k\to 0$ as $k\to \infty$
needs to be refined when $D$ has (small) positive eigenvalues for which
$\gamma_i\approx 0$.
  Components $i$ for which $\gamma_i\approx 0$
  will be called
\\
  \centerline{\em critical components}
\\
  below. 
  Given the choice
  $\beta = (1-\sqrt{\underline \alpha})^2$,
  critical components are those with $D_{i,i}$ close to
  $\underline m$.

%\newpage
\subsection{The Schur canonical form}\label{sec.schur}

It turns out that the case $\gamma_i\approx 0$
can be analyzed by using the Schur-decomposition
of $M^{(i)}$ with an upper right triangular matrix $R$
and a transformation matrix $T$ given by
\begin{align*}
 TRT^{-1}:=\begin{pmatrix}1&0\\ \lambda_+&1\end{pmatrix}\begin{pmatrix}\lambda_+&1\\ 0&\lambda_-\end{pmatrix}\begin{pmatrix}1&0\\ -\lambda_+&1\end{pmatrix}=\begin{pmatrix}0&1\\ -\lambda_+\lambda_-&\lambda_++\lambda_-\end{pmatrix}=\begin{pmatrix}0&1\\ -\beta&\beta_i\end{pmatrix}=M^{(i)}.
\end{align*}
Note that this decomposition is valid for both cases: for $\gamma_i=0$
as well as for $\gamma_i\neq 0$. 
(When $\gamma_i=0$ i.e.\ when $\lambda_+=\lambda_-$
this is the Jordan canonical form.)

Now, in place of \mr{argu1}, also the argument 
\bee\label{argu2}
  \|(M^{(i)})^{k}\| = \|(TRT^{-1})^{k}\|
  = \|TR^{k}T^{-1}\| 
 \leq \|T\| \cdot \|T^{-1}\| \cdot \|R^{k}\|  =
\cond(T)\cdot \|R^{k}\| 
\ene
can be used. While both \mr{argu1} and \mr{argu2} lead to valid bounds
for $\|(M^{(i)})^{k}\|$ whenever \mr{argu1} is defined and thus,
the better of both bounds can be used, the estimate below relies on
\mr{argu2} which is always well defined.
Bounds on $\cond(T)$ and on $\|R^{k}\|$ are derived next:

By Gershgorin's theorem, given a triangular matrix
$$
\hat R =
\begin{pmatrix}a & b\\ 0 & c\end{pmatrix}
$$
with complex numbers $a,b,c$ and $|a|\ge |c|$, the 2-norm of
$\hat R$ is bounded by
\bee\label{2norm}
\|\hat R\|_2 =
\lambda_{max}\left(
\begin{pmatrix}\bar a & 0\\ \bar b & \bar c\end{pmatrix}
\begin{pmatrix}a & b\\ 0 & c\end{pmatrix}
\right)^{1/2}
\!\!=
\lambda_{max}\left(
\begin{pmatrix}|a|^2 & \bar ab\\ a\bar b & |c|^2+|b|^2\end{pmatrix}
\right)^{1/2}
\le 
\sqrt{|a|^2+|ab|+|b|^2}
\ene
(with $\lambda_{max}$ denoting the maximum eigenvalue).
%Gershgorin's theorem is not the same for $\hat R$ and $\hat R^H$ but 

Since $\|\hat R\|_2=\|\hat R^H\|_2$ the bound \mr{2norm} also applies to
the matrix $T$ and yields
\bee\label{est.t}
\|T\|_2^2
\le 1+|\lambda_+|+|\lambda_+|^2 \le 3
\ene
since $|\lambda_+|\le \rho(M^{(i)})< 1$.
The same estimate applies to $\|T^{-1}\|_2^2$, so that
$\cond(T)\le 3$.

Inductively, it can be verified that the $k^\text{th}$ power of the
matrix $R$ is given by
\begin{align*}
 R^k = \begin{pmatrix} \lambda_+^k & \sum_{\ell=0}^{k-1}\lambda_+^\ell\lambda_-^{k-1-\ell} \\ 0 & \lambda_-^k \end{pmatrix}.
\end{align*}
%Without loss of generality assume that $|\lambda_-|\le|\lambda_+|$, else rename $\lambda_-$ and $\lambda_+$.
By definition $|\lambda_+|\leq\rho$ and $|\lambda_-|\leq\rho$ apply. (If $\gamma_i=0$, then actually $\rho=|\lambda_-|=|\lambda_+|$ ---
and this is the ``worst case'' in the estimate below).
In any case, it follows with the triangle inequality
that the absolute value of the upper right entry
of $R^k$ is bounded by 
$|(R^k)_{1,2}|\le k\rho^{k-1}$, and then, 
by \mr{2norm}, that
\bee\label{est.r}
\|R^k\|\le \rho^{k-1}\sqrt{ \rho^2+ k\rho+k^2}
\le \rho^{k-1}\sqrt{ 1+ k+k^2}\le \rho^{k-1}(k+1).
\ene
Since the condition number of $T$ is at most 3, the norm 
$\|(M^{(i)})^{k}\|$ is at least $\frac 13 \|R^k\|$ and at most $3\|R^k\|$.
The latter observation can be used to bound the
number of iterations needed to reduce the (unknown) Euclidean distance
of a given point $x^0$ to the optimal solution by a factor $\epsilon>0$.
This is done in the next section with a slightly tighter bound
of $\|(M^{(i)})^k\|$.

\subsection{On the rate of convergence}
Inserting \mr{est.t} and \mr{est.r} in \mr{argu2} yields
$
\|(M^{(i)})^k\|_2\le 3 \rho^{k-1}(k+1)
$.
%where $\rho=\rho(M^{(i)})\ge |\lambda_+|$. 
This bound is based on
the sub-multiplicativity of the 2-norm. 
When computing $(M^{(i)})^k$ explicitly, it can be reduced to
\bee\label{verylast}
\|(M^{(i)})^k\|_2\le 2 \rho^{k-1}(k+1).
\ene
The somewhat tedious calculations leading to \mr{verylast} are
given in the Appendix.
Let a desired accuracy $\epsilon\in(0,1)$ be given.
To obtain a sufficient condition for
\bee\label{boundk}
2 \rho^{k-1}(k+1)\le \epsilon  \ \ \Longleftrightarrow \ \
(k\!-\!1)\ln(\rho)  +\ln(k+1)\le \ln(\tfrac \epsilon 2)
\ene
the bound $\ln(\rho)\le \rho-1$ is used,
\beea
\mr{boundk}   \ \ &\Longleftarrow& \ \
(k\!-\!1)(\rho-1)  +\ln(k+1)\le \ln(\tfrac \epsilon 2) \nonumber \\
&\Longleftrightarrow& \ \ (1-\rho)(k\!-\!1) \ge \ln(\tfrac 2\epsilon)
+ \ln(k+1)\label{bound1}
\enea
Denote
$$
\delta := \tfrac1{\sqrt{2\ \ocond(D)}} =
 \tfrac{1-\sqrt\beta}{2} \le \tfrac{1-\rho}{2}.
$$
If the following two
inequalities are satisfied then \mr{bound1} is satisfied as well,
\bee\label{bound2}
\delta (k-1) \ge \ln(\tfrac 2\epsilon) \qquad \hbox{and} \qquad
\delta (k-1) \ge \ln(k+1).
\ene
Here $\delta>0$ depends on $\beta=(1-\sqrt{\underline{\alpha}})^2$,
i.e.\ on the choice of the parameters $\alpha,\beta$ that in turn
depend on the bounds $\underline{m}$ and $\overline{m}$
of the eigenvalues of $D$. Below it is assumed that $\delta \le \exp(-2)$
-- which is the case when $\ocond(D)\ge 28$.

Then, the second relation in \mr{bound2} is satisfied whenever
$k\ge \bar k:=\frac 2\delta \ln(\frac 1\delta)-1$.

\noindent
Indeed,
\beeas
\delta(\bar k-1) - \ln(\bar k + 1)&=&
2\ln(\tfrac 1\delta)-2\delta-\ln(\tfrac 2\delta \ln(\tfrac 1\delta))\\
&=& \ln(\tfrac 1\delta)-2\delta-\ln(\ln(\tfrac 1\delta))-\ln(2)\\
&>&0
\eneas
for $0<\delta \le \exp(-2)$.

Thus it suffices to ensure that the first inequality in \mr{bound2}
implies the second. This is the case when $\epsilon$ is sufficiently small:
Indeed, the first bound in \mr{bound2} is satisfied
for $k=\bar k$ with equality if 
$$
\epsilon = \bar\epsilon :=
2\delta^2e^{2\delta} .
$$
For $\epsilon \le \bar\epsilon$ all values of $k$ satisfying
the first relation in \mr{bound2} are greater or equal to $\bar k$
so that the second relation is satisfied as well.
Since $e^{2\delta}>1$ it is sufficient to choose
\blue{$\epsilon \le 2\delta^2 = 1\, / \ocond(D)$}. The first relation in
\mr{bound2} can be written as 
\bee\label{last}
k\ge \ln\left(\tfrac 2\epsilon \right)\cdot \tfrac 1\delta +1 = \ln\left(\tfrac 2\epsilon \right)\cdot
\sqrt{ 2\ \ocond(D)}+1.
\ene 
Thus, when $k$ satisfies \mr{last}
with $\epsilon \le 2\delta^2$, relation \mr{boundk} is satisfied and
\blue{
\beeas
&&\sqrt 2 \| \tfrac 12 (x^{k-1}+x^k)\|_2
=\left\| \tfrac 12\left(\left(
\begin{matrix} x^{k-1} \\ x^k \end{matrix} \right)+\left(
\begin{matrix} x^k \\ x^{k-1} \end{matrix} \right)\right)
\right\|_2\le
\tfrac 12\left(\left\| \left(\begin{matrix} x^{k-1} \\ x^k \end{matrix} \right)
\right\|_2 + \left\| \left(\begin{matrix} x^{k} \\ x^{k-1} \end{matrix} \right)
\right\|_2\right) \\
&=&  
\tfrac 12\left(\left\| \left(\begin{matrix} x^{k-1} \\ x^k \end{matrix} \right)
\right\|_2 + \left\| \left(\begin{matrix} x^{k-1} \\ x^k \end{matrix} \right)
\right\|_2\right) =
\left\| \left(\begin{matrix} x^{k-1} \\ x^k \end{matrix} \right)
\right\|_2 \le \epsilon
\left\| \left(\begin{matrix} x^{0} \\ x^1 \end{matrix} \right)
\right\|_2 \le \sqrt 2 \epsilon \|x^0\|_2
\eneas
}%
where the last inequality follows from $\|x^1\|_2\le \|x^0\|_2$,
the first step being a plain steepest descent step with step length
$2/\overline m$.
Summarizing, the claim of Theorem \ref{thm1} follows.
\hfill $\square$

\blue{
\subsection{Nesterov's accelerated gradient method}\label{sec.nest}
For the analysis of Nesterov's accelerated gradient method
again it suffices to consider the function $f$ in \mr{assu}
with $\nabla f(x) = Dx$.
Replacing $(HBM)$ with $(NAG)$ in the form \mr{nag} it follows that the matrix
$
M^{(i)}=\left(\begin{matrix} 0 && 1 \\ -\beta  && \beta_i \end{matrix} \right)
$
in \mr{defmi} %for $(HBM)$ 
is to be replaced with 
$$
M^{(i)}:=\left(\begin{matrix} 0 && 1 \\ -\beta(1-\alpha_i)  && (1+\beta)(1-\alpha_i) \end{matrix} \right)
$$
where again, $\alpha_i:=\alpha D_{i,i}$.
Thus, replacing $\beta$ and $\beta_i$ in \mr{defmi} with
$\beta(1-\alpha_i)$ and $(1+\beta)(1-\alpha_i)$ it follows that $\rho$
in \mr{defrho} changes to
\bee\label{defrho2}
\rho:=\rho(M^{(i)}) = \left\{\begin{array}{ll}
\sqrt{\beta(1-\alpha_i)} & \hbox{if} \ \ (1+\beta)^2(1-\alpha_i)^2\le 4\beta(1-\alpha_i)\\[4pt]
\frac 12\left( | (1+\beta)(1-\alpha_i) |+\bar\gamma_i \right)
& \hbox{else}
\end{array}\right.
\ene
where $\bar \gamma_i:= \sqrt{(1+\beta)^2(1-\alpha_i)^2- 4\beta(1-\alpha_i)}$.
The first case of definition \mr{defrho2} can only be obtained for
$\alpha_i\le 1$. 
This results in
half the step length $\alpha:=\frac 1{\overline{m}}$
compared to \mr{alphapol} for $(HBM)$ and 
it then follows for $\alpha_i\in [1/\ocond(D),1]$ and
$\beta:=\tfrac{\left(\sqrt{\ocond(D)}-1\right)^2}{\ocond(D)-1}$ that
$(1+\beta)^2(1-\alpha_i)\le 4\beta$
and
$$
\rho=\sqrt{\beta(1-\alpha_i)}\le 1-\tfrac 1{\sqrt{\ocond(D)}}.
$$
This is the same bound as for the $(MM)$ when taking half the step length
$\alpha$ or -- which is the same -- when replacing $\overline{m}$
in $(MM)$ with $2\overline{m}$. 
In contrast to the $(MM)$, however, 
based on this analysis a theoretical justification of a longer step
length $\alpha > 1/\overline{m}$ is not possible.\\
The eigenvalues of $M^{(i)}$ are now given by
$$
\lambda_\pm = \tfrac 12 (\ (1+\beta)(1-\alpha_i)\pm
\sqrt{(1+\beta)^2(1-\alpha_i)^2-4\beta(1-\alpha_i)}\ )
$$
and with this definition 
the matrix $T$ transforming $M^{(i)}$ to  the Schur canonical form
is the same as in Section \ref{sec.schur}.
The iteration complexity for $(MM)$ thus applies to $(NAG)$
as well when replacing $\ocond(D)$ with $2\,\ocond(D)$.
Summarizing, the following result is obtained:
\begin{theorem}\label{thm2}
  Let $f:\re^n\to\re$ be a strictly convex quadratic function
  with minimizer $x^*$ and
  denote $H:=\nabla^2f(x^*)$. Let
  $\overline{m}$ be an upper bound for the eigenvalues of $H$
  and let $0<\underline{m}$ be a lower bound for the eigenvalues
  of $H$. Then $\ocond(H):=\overline{m}/\underline{m}$
  is an upper bound for the condition number of $H$. Assume
  that $\ocond(H)\ge 28$.
  \\
  In $(HBM)$ define $\alpha=1/\overline{m}$, \
  $\underline{\alpha}=1/\ocond(H)$,  \ 
  $\beta = \tfrac{\left(1-\sqrt{\underline{\alpha}}\right)^2}
  {1-\underline{\alpha}}$
  and let  
  $\epsilon\le 1\, / \, \ocond(H)$ be given.
  Then, given $x^0\in\re^n$, after
  $$
  1+\lceil 2\,\sqrt{\ocond(H)}\
  \ln(\tfrac {2}\epsilon)\rceil
  $$
  steps of the $(NAG)$, an approximate solution
  $\bar x^k:=\tfrac 12 (x^{k-1}+x^k)$ is
  generated with
  $$
  \|\bar x^k-x^*\|_2\le
  %\sqrt{\|x^k-x^*\|_2^2+ \|x^{k-1}-x^*\|_2^2} \le
  \epsilon \|x^0-x^*\|_2.
  $$
\end{theorem}
}

\noindent
\subsection{A note on the rate of convergence}
Theorem~\ref{thm1} and Theorem~\ref{thm2} address the iteration
complexity with respect to the Euclidean norm while the
asymptotic rate of convergence of $(NAG)$ or $(HBM)$ (with respect to any
fixed norm)
is given by the maximum spectral radius of all matrices $M^{(i)}$ which coincides
with the value $1-\sqrt{\underline\alpha}$ under the assumptions of
Theorems~\ref{thm1} and \ref{thm2}.
Thus, the asymptotic rate of convergence is
\begin{itemize}
\item
$1- \sqrt{1/\ocond(H)}$ for $(NAG)$ and
\item
$1- \sqrt{2/\ocond(H)}$ for $(HBM)$.
\end{itemize}
Here, the rate for $(HBM)$ is somewhat worse than the asymptotic rate
\hbox{$1-2/\sqrt{\ocond(H)}$} 
for $(HBM)$ established in \cite{polyak}, Theorem 9 (3),
the difference being due to the limitation of the step length $\alpha$
which is nearly twice as long in \cite{polyak} as considered here
-- and four times as long as for Nesterov's accelerated gradient method.

The rates of convergence are not just upper bounds but exact rates
once the values $\underline m$ and $\overline m$ are given,
and hence, 
the difference in the rates of convergence for $(NAG)$ and $(HBM)$
can be explained just by the limitations of the step length: when the $(HBM)$
is restricted to steps $\alpha = 1/\overline m$ as for $(NAG)$, its asymptotic
rate of convergence is exactly the same as $(NAG)$.

The above rate for Nesterov's accelerated gradient method is the same
as established for the convergence of the function values in
\cite{nestg2}, Theorem 2.2.3.
Translated to the error measure
$\|x-x^*\|_H = \sqrt{f(x)-f(x^*)}$ the rate of convergence
of \cite{nestg2}, Theorem 2.2.3 reduces to approximately
$1-\sqrt{1/2\, \ocond(H)}$.

\section{Conclusion}
The asymptotic rate of convergence of the momentum method for
fixed parameters has been analyzed in the classical work \cite{polyak}.
Here, the constants in this work are considered and 
an explicit bound for the iteration complexity is derived
\blue{that also applies in slightly modified form to
Nesterov's accelerated gradient method}.

\subsection*{Acknowledgment}
The authors would like to thank Robert Gower for helpful
e-mail-discussions.

\subsection*{Appendix}
The matrix $M^{i,k}:=(M^{(i)})^k$ can also be written as
\begin{align*}
(M^{(i)})^k &= TR^kT^{-1} =
\begin{pmatrix}1&0\\\lambda_+&1\end{pmatrix}
\begin{pmatrix}\lambda_+^k&\ \ \ \sum_{\ell=0}^{k-1}\lambda_+^\ell\lambda_-^{k-1-\ell}\\[4pt]
  0&\lambda_-^k\end{pmatrix}
\begin{pmatrix}1&0\\-\lambda_+&1\end{pmatrix}\\
&=
\begin{pmatrix}1&0\\\lambda_+&1\end{pmatrix}
\begin{pmatrix}\lambda_+^k-\lambda_+\sum_{\ell=0}^{k-1}\lambda_+^\ell\lambda_-^{k-1-\ell}&
  \ \ \ \ \ \sum_{\ell=0}^{k-1}\lambda_+^\ell\lambda_-^{k-1-\ell}\\[4pt]
  -\lambda_+\lambda_-^k&\lambda_-^k\end{pmatrix}\\
&=
\begin{pmatrix}\lambda_+^k-\lambda_+\sum_{\ell=0}^{k-1}\lambda_+^\ell\lambda_-^{k-1-\ell}&
  \sum_{\ell=0}^{k-1}\lambda_+^\ell\lambda_-^{k-1-\ell}\\[4pt]
  \lambda_+^{k+1}-\lambda_+^2\sum_{\ell=0}^{k-1}\lambda_+^\ell\lambda_-^{k-1-\ell}-\lambda_+\lambda_-^k&
  \ \ \ \ \ \lambda_+\sum_{\ell=0}^{k-1}\lambda_+^\ell\lambda_-^{k-1-\ell}+\lambda_-^k\end{pmatrix}\\
&=
\begin{pmatrix}\lambda_+^k-\sum_{\ell=0}^{k-1}\lambda_+^{\ell+1}\lambda_-^{k-1-\ell}&
  \sum_{\ell=0}^{k-1}\lambda_+^\ell\lambda_-^{k-1-\ell}\\[4pt]
  \lambda_+^{k+1}-\sum_{\ell=0}^{k-1}\lambda_+^{\ell+2}\lambda_-^{k-1-\ell}-\lambda_+\lambda_-^k&
  \ \ \ \ \ \sum_{\ell=0}^{k-1}\lambda_+^{\ell+1}\lambda_-^{k-1-\ell}+\lambda_-^k\end{pmatrix}\\
&=
\begin{pmatrix}-\sum_{\ell=0}^{k-2}\lambda_+^{\ell+1}\lambda_-^{k-1-\ell}&
  \sum_{\ell=0}^{k-1}\lambda_+^\ell\lambda_-^{k-1-\ell}\\[4pt]
  -\sum_{\ell=-1}^{k-2}\lambda_+^{\ell+2}\lambda_-^{k-1-\ell}&
  \ \ \ \sum_{\ell=-1}^{k-1}\lambda_+^{\ell+1}\lambda_-^{k-1-\ell}\end{pmatrix}\\
&=
\begin{pmatrix}-\sum_{t=1}^{k-1}\lambda_+^{t}\lambda_-^{k-t}&
  \ \ \sum_{\ell=0}^{k-1}\lambda_+^\ell\lambda_-^{k-1-\ell}\\[4pt]
  -\sum_{t=0}^{k-1}\lambda_+^{t+1}\lambda_-^{k-t}&
  \sum_{t=0}^{k}\lambda_+^{t}\lambda_-^{k-t}\end{pmatrix}
  =: \begin{pmatrix}p&q\\r&s\end{pmatrix}\in\mathbb{R}^{2\times2},
\end{align*}
where in the last row the index shift $t:=\ell+1$ was used.
(The definition of the real numbers $p,q,r,s$ depends on the
  possibly complex eigenvalues $\lambda_+$ and $\lambda_-$.)

The matrix product $(M^{i,k})^\text{T}M^{i,k}$ is
\begin{align*}
 (M^{i,k})^\text{T}M^{i,k}
 = \begin{pmatrix}p&r\\q&s\end{pmatrix}\begin{pmatrix}p&q\\r&s\end{pmatrix}
 = \begin{pmatrix}p^2+r^2&pq+rs\\pq+rs&q^2+s^2\end{pmatrix}
 = \begin{pmatrix}a&b\\b&d\end{pmatrix}
 \end{align*}
 with
 \begin{align*}
  a&:=p^2+r^2 = \left(\sum_{\ell=1}^{k-1}\lambda_+^{\ell}\lambda_-^{k-\ell}\right)^2+\left(\sum_{\ell=0}^{k-1}\lambda_+^{\ell+1}\lambda_-^{k-\ell}\right)^2,\\
  b&:= pq+rs = -\left(\sum_{\ell=1}^{k-1}\lambda_+^{\ell}\lambda_-^{k-\ell}\right)\left(\sum_{\ell=0}^{k-1}\lambda_+^\ell\lambda_-^{k-1-\ell}\right)
  -\left(\sum_{\ell=0}^{k-1}\lambda_+^{\ell+1}\lambda_-^{k-\ell}\right)\left(\sum_{\ell=0}^{k}\lambda_+^{\ell}\lambda_-^{k-\ell}\right),\\
  d&:= q^2+s^2 = \left(\sum_{\ell=0}^{k-1}\lambda_+^\ell\lambda_-^{k-1-\ell}\right)^2+\left(\sum_{\ell=0}^{k}\lambda_+^{\ell}\lambda_-^{k-\ell}\right)^2.
 \end{align*}
 By definition of the spectral radius $\rho=\rho(M^{(i)})<1$ the
 inequalities $|\lambda_+|\leq\rho$ and $|\lambda_-|\leq\rho$ hold, and thus,
$$
a \le 2k^2\rho^{2k}\quad  \hbox{and} \quad d \le 2(k+1)^2\rho^{2k-2}. 
$$
Since $(M^{i,k})^\text{T}M^{i,k}\succeq 0$ 
the maximum eigenvalue of $(M^{i,k})^\text{T}M^{i,k}$
is at most equal to the trace, i.e. 
$$
\lambda_{\max}((M^{i,k})^\text{T}M^{i,k}) 
\le a+d \le 4\rho^{2(k-1)} (k+1)^2.
$$
This way, an upper bound for the norm of $M^{i,k}$ is given by
\bee\label{newbdmik}
 \|M^{i,k}\|_2 \leq  2\rho^{k-1} (k+1).
\ene
To estimate in how far the above bound is tight, consider the
case of $\gamma_i=0$, i.\,e. $\rho=|\lambda_+|=|\lambda_-|$.
It then follows
$$
a \ge 2(k-1)^2\rho^{2k+2},\qquad  b \le -2(k-1)^2\rho^{2k+1},\qquad
d \ge 2(k-1)^2\rho^{2k} 
$$
and with $z:=(1,\ -1)^T$ that
$$
\lambda_{\max}((M^{i,k})^\text{T}M^{i,k})\ge
\tfrac{z^T(M^{i,k})^\text{T}M^{i,k}z}{z^Tz}\ge \tfrac{8(k-1)^2\rho^{2k+2}}{2}.
$$
Thus,
$$
\|M^{i,k}\|_2 \geq  2\rho^{k+1} (k-1),
$$
i.e.\ up to the changes  $k+1\longleftrightarrow k-1$
the bound \mr{newbdmik} is sharp.

Again, the ``critical components''
  i.e.\ those with $\gamma_i\approx 0$ lead to a poor convergence,
  i.e. to the case where \mr{newbdmik} is almost sharp.
  On the other hand, $(M^{i,k})^TM^{i,k}$ is almost singular, and 
  for a critical component $i$ the quantity $\alpha_i$
  is small so that the iterates $x^0$ and $x^1$ of $(MM)$ with
  step length $2/\overline m$ satisfy $x^0_i \approx x^1_i$. 
  This again implies that
  the vector $z_{(i)}^1:=(x^0_i,\ x^1_i)^T$ is close to the eigenvector
  of $(M^{i,k})^TM^{i,k}$ associated with the small eigenvalue, i.e.
  $$
  \|M^{i,k}z_{(i)}^1\|_2^2=(z_{(i)}^1)^T(M^{i,k})^TM^{i,k}z_{(i)}^1
  \ll \|M^{i,k}\|_2^2 \|z_{(i)}^1\|_2^2.
  $$

%\newpage


\begin{thebibliography}{99}

\bibitem{barre}
Barr\'e, M., Taylor, A., d'Aspremont, A. (2020):
Complexity Guarantees for Polyak Steps with Momentum.
33rd Annual Conference on Learning Theory,
Proceedings of Machine Learning Research, Vol. 125, 1-27.

\bibitem{defazio_mm}
Defazio, A. (2021):
Momentum via Primal Averaging: Theoretical
Insights and Learning Rate Schedules for
Non-Convex Optimization.
https://arxiv.org/pdf/2010.00406.pdf

\bibitem{defossez}
%Alexandre Défossez defossez@meta.com Meta AI
%Léon Bottou Meta AI Francis Bach INRIA / PSL Nicolas Usunier Meta AI
D\'efossez, A., Bottou, L., Bach, F., Usunier, N. (2022):
A Simple Convergence Proof of Adam and Adagrad.
Transactions on Machine Learning Research.
https://arxiv.org/pdf/2003.02395.pdf

\bibitem{diakonikolas}
Diakonioklas, J., Jordan, M.I. (2021):
Generalized Momentum-Based Methods: a Hamiltonian Perspective.
SIAM J.Optim. Vol. 31, No. 1, 915-944.
%JELENA DIAKONIKOLAS † AND MICHAEL I. JORDAN ‡

\bibitem{donoghue}
O’Donoghue, B.,  Cand\`es. E. (2015):
Adaptive Restart for Accelerated Gradient Schemes.
Foundations of computational mathematics,
Vol. 15, No 3, 715-732.

\bibitem{ganesh}
Ganesh, S., Deb, R., Thoppe, G., Budhiraja, A. (2022):
Does Momentum Help in Stochastic Optimization?
A Sample Complexity Analysis.
https://arxiv.org/abs/2110.15547v3

\bibitem{gitman}
%Igor Gitman Hunter Lang Pengchuan Zhang Lin Xiao
Gitman, I., Lang, H., Zhang, P., Xiao, L. (2019):
Understanding the Role of Momentum in Stochastic Gradient Methods.
H. Wallach, H. Larochelle, A. Beygelzimer, F. d'Alch\'e-Buc, E. Fox,
R. Garnett (eds):
Advances in Neural Information Processing Systems, 32 (NeurIPS 2019),
https://proceedings.neurips.cc/paper/2019

\bibitem{goh}
Goh, G. (2017):
Why Momentum Really Works.
Distill,
http://distill.pub/2017/momentum

%\bibitem{goodf}
%Ian Goodfellow, Yoshua Bengio, and Aaron Courville. Deep learning. MIT press, 2016.  
%Goodfellow, I., Bengio, Y., Courville, A. (2016):
%Deep learning. MIT press.

\bibitem{gratton}
Gratton, S., Jerad, S.,  Toint, Ph. L. (2022):
First-Order Objective-Function-Free Optimization Algorithms
and Their Complexity.
  
\bibitem{gratton2}
Gratton, S., Jerad, S.,  Toint, Ph. L. (2022):
Parametric complexity analysis for
a class of first-order Adagrad-like algorithms.
https://arxiv.org/pdf/2203.01647.pdf

\bibitem{gratton3}
Gratton, S.,  Toint, Ph. L. (2022):
OFFO minimization algorithms for second-order optimality and their complexity.
https://arxiv.org/pdf/2203.03351.pdf

%\bibitem{hagedorn1}
%Hagedorn, M., Jarre, F. (2022):
%Optimized convergence of stochastic gradient descent by weighted averaging.
%https://arxiv.org/abs/2209.14092

\bibitem{li}
Li, R.C.  (2008):
On Meinardus' examples for the conjugate gradient method,
Mathematics of Computation,
Vol. 77, Nr. 261,  335-352.

\bibitem{nestg}
Nesterov, Y.E. (1983):
A method for solving the convex programming problem with convergence rate $O(1/k^2)$.
Dokl. akad. nauk Sssr 269, 543-547.

\bibitem{nestg2}
Nesterov, Y.E. (2003):
Introductory lectures on convex optimization: A basic course.
Springer Science and Business Media,
Vol. 87.

\bibitem{polyak}
Polyak, B.T. (1964):
Some methods of speeding up the convergence of iteration methods.
USSR Computational Mathematics and Mathematical Physics, 4(5):1-17.
  
\bibitem{gower1}
Sebbouh, O., Gower, R.M., Defazio, A. (2020):
On the convergence of the Stochastic Heavy Ball Method.
https://othmanesebbouh.github.io/publications/heavy$\underline{\ }$ball.pdf

\bibitem{hagedorn2}
The GitHub repository: https://github.com/MHagedorn/momentum (2022)

\bibitem{yang}
%Tianbao Yang Qihang Lin Zhe Li
Yang, T., Lin, Q., Li, Z. (2016):
Unified Convergence Analysis of Stochastic Momentum
Methods for Convex and Non-convex Optimization. 
https://arxiv.org/pdf/1604.03257.pdf



\end{thebibliography}
\end{document}